\documentclass[10pt]{amsart}
\usepackage{amsmath,amssymb,amscd}

\usepackage{enumerate}
\usepackage[all]{xy}
\usepackage{amssymb}
\usepackage{mathrsfs}
\usepackage{stmaryrd}
\usepackage{enumerate}
\usepackage{epsfig}
\usepackage{layout}
\newtheorem{thm}{Theorem}
\newtheorem{theoremalpha}{Theorem}

\newtheorem{questionalpha}{Question}

\newtheorem{question*}[]{Question}

\theoremstyle{definition}

\newtheorem{defi}[thm]{Definition}

\newenvironment{sis}{\left\{\begin{aligned}}{\end{aligned}\right.}

\def\w{\widetilde}
\def\mbb{\mathbb}
\def\mcl{\mathcal}

\def\ov{\overline}

\def\inj{\hookrightarrow}
\def\surj{\twoheadrightarrow}

\def\Z{\mbb Z}

\def\cO{\mcl O}

\def\deg{{\rm deg}}

\def\Jac{\rm{Jac}}

\def\P{\mathbb{P}}

\def\Mgb{\bar{M}_g}
\def\Mgps{ \overline{M}^{\rm{ps}}_g}

\def \Hilb{{\rm Hilb}}
\def \Chow{{\rm Chow}}
\def \Ch{{\rm Ch}}

\def \Gr{{\rm Gr}}
\def\SL{\rm{SL}}

\def\Aut{{\rm Aut}}

\def \OO{\mcl O}
\def \Sym{{\rm Sym}}

\def \exc{X_{\rm exc}}

\def \Jst{\ov{J}_{d,g}}

\def \JJps{\ov{\mathcal J}_{d,g}^{\rm ps}}
\def \Jps{\ov{J}_{d,g}^{\rm ps}}

\def \Pdgps{\ov{J}_{d,g}^{\rm ps}}

\def\bar{\overline}

\input xy
\xyoption{all}

\begin{document}

\normalsize

%\vskip 1in

\title
%[GIT for polarized curves]
[On GIT quotients of Hilbert and Chow schemes of curves]
{On GIT quotients of Hilbert and Chow schemes of curves}

\author{Gilberto Bini, Margarida Melo, Filippo Viviani}

\email{gilberto.bini@unimi.it} \curraddr{{\sc Dipartimento di Matematica \\ Universit\`a degli Studi di Milano \\ Via C. Saldini 50 \\ 20133 Milano \\ Italy.}}

\email{mmelo@mat.uc.pt}\curraddr{{\sc Departamento de Matem\'atica \\ Universidade de Coimbra\\
Largo D. Dinis, Apartado 3008 \\ 3001 Coimbra \\Portugal.}}

\email{viviani@mat.uniroma3.it} \curraddr{{\sc Dipartimento di Matematica \\ Universit\`a Roma Tre \\ Largo S. Leonardo Murialdo 1 \\ 00146 Roma \\ Italy.} }
%and {\sc Departamento de Matem\'atica \\ Universidade de Coimbra\\
%Largo D. Dinis, Apartado 3008 \\ 3001 Coimbra \\Portugal.}}

%\thanks{The first named author has been partially supported by "FIRST" Universit\`{a} di Milano and by MIUR Cofin 2008 %- Variet\`{a} algebriche: geometria, aritmetica, strutture di Hodge. The second named author has been partially %supported by GNSAGA of INdAM and by MIUR Cofin 2008 - Geometria delle variet\`{a} algebriche e dei loro spazi di %moduli. The third named author has been partially supported by FCT-Ci\^encia2008 and by MIUR Cofin 2008 - Geometria %delle variet\`{a} algebriche e dei loro spazi di moduli.}

\date{\today}

\begin{abstract}
The aim of this note is to announce some results on the GIT problem for the Hilbert and Chow scheme of curves of degree $d$ and genus $g$ in the projective 
space of dimension $d-g$, whose full details will appear in \cite{BMV}.
%as their degree $d$ decreases with respect to their genus $g$.
% In particular, we aim at a description of the semistable, polystable, stable points and how semistable orbits get identified in the quotient.
In particular, we extend the previous results of L. Caporaso up to $d>4(2g-2)$ and we observe that this is sharp.
In the range $2(2g-2)<d<\frac{7}{2} (2g-2)$, we get a complete new description of the GIT quotient.
As a corollary, we get a new compactification of the universal Jacobian   over the moduli space of pseudo-stable curves.
\end{abstract}

\thanks{The first named author has been partially supported by ``FIRST'' Universit\`a di Milano and by MIUR PRIN 2010 - Variet\`a algebriche: geometria, aritmetica e strutture di Hodge.
The second named author was partially supported by CMUC and FCT (Portugal) through European program COMPETE/FEDER, 
by the FCT project \textit{Espa\c cos de Moduli em Geometria Alg\'ebrica} (PTDC/MAT/111332/2009), by the FCT project \textit{Geometria Alg\'ebrica em Portugal} (PTDC/MAT/099275/2008) and by the Funda\c c\~ao Calouste Gulbenkian program ``Est\'imulo \`a investiga\c c\~ao 2010''.
The third named author was partially supported by CMUC and FCT (Portugal) through European program COMPETE/FEDER and 
by the FCT project \textit{Espa\c cos de Moduli em Geometria Alg\'ebrica} (PTDC/MAT/111332/2009).}
%and by MIUR Cofin 2008 - \textit{Geometria delle variet\`{a} algebriche e dei loro spazi di moduli}.}

\maketitle

\section{Motivation}

One of the first successful applications of Geometric Invariant Theory (GIT for short) was the construction of the moduli space $M_g$ of smooth curves of genus $g\geq 2$ together with its compactification $\Mgb$ via stable curves,
carried out by Mumford (\cite{Mum}) and Gieseker (\cite{Gie}).
Indeed, the moduli space of stable curves was constructed as a GIT quotient of the locally closed subset of a suitable Hilbert scheme
(as in \cite{Gie}) or  Chow scheme (as in \cite{Mum})
parametrizing $n$-canonically embedded curves, for $n\geq 5$ (see also \cite[Chap. 4, Sec. C]{HM} or \cite[Sec. 3]{Mo}).

Recently, there has been a lot of interest in extending the above GIT analysis to smaller values of $n$, specially in connection with the so called Hassett-Keel program whose ultimate goal
is to find the minimal model of $M_g$ via the successive constructions of modular birational models of $\Mgb$ (see \cite{FS} and \cite{AH} for nice overviews).
The first work in this direction is due to Schubert,
who described in \cite{Sch}  the GIT quotient of the locus of $3$-canonically embedded curves (of genus $g\geq 3$) in the Chow scheme
as the coarse moduli space $\Mgps$ of pseudo-stable curves (or \emph{p-stable} curves for short), i.e., reduced, connected, projective curves with finite
automorphism group, whose only singularities are nodes and ordinary cusps, and which have no elliptic tails. Indeed, it was shown by Hyeon-Morrison in \cite{HMo} that
one gets the same GIT quotient for the Hilbert scheme of $3$-canonically embedded curves and for the Hilbert or Chow scheme of $4$-canonically embedded curves.
Later, Hasset-Hyeon constructed in \cite{HH1} a modular map $T:\Mgb\to \Mgps$ which on geometric points sends a stable curve onto the p-stable curve obtained by contracting all its elliptic
tails to cusps. Moreover the authors of loc. cit. identified the map $T$ with the first contraction in the Hassett-Keel program for $\Mgb$.
Finally the GIT quotient of the Hilbert and Chow scheme of $2$-canonically embedded curves was studied in great detail by Hassett-Hyeon in \cite{HH2}. For some partial results on the GIT
quotient for the Hilbert scheme of $1$-canonically  embedded curves, see the works of Alper-Fedorchuk-Smyth (\cite{AFS}, \cite{AFS2}) and of Fedorchuk-Jensen (\cite{FJ}).

From the point of view of constructing new projective birational models of $\Mgb$, it is of course natural to restrict the GIT analysis to the locally closed subset inside the Hilbert or the Chow
scheme parametrizing $n$-canonical embedded curves. However, the problem of describing the whole GIT quotient seems very natural and interesting too. The first result in this direction is the
breakthrough work of Lucia Caporaso \cite{Cap}, where she describes the GIT quotient of the Hilbert scheme of connected curves of genus $g\geq 3$ and degree $d\geq 10(2g-2)$
in $\P^{d-g}$.
The GIT quotient that she obtains is indeed a modular compactification $\ov{J}_{d,g}$ of the universal Jacobian $J_{d,g}$, which is the moduli scheme parametrizing pairs $(C,L)$ where
$C$ is a smooth curve of genus $g$ and $L$ is a line bundle on $C$ of degree $d$. Note that recently Li-Wang in \cite{LW} have given a different proof of the Caporaso's result for $d\gg0$.

Our work is motivated by the following

\vspace{0,2cm}

\textbf{Problem(I):}
\emph{Describe the GIT quotient of the Hilbert and Chow scheme of curves of genus $g$ and degree $d$ in $\P^{d-g}$, as $d$ decreases with respect to $g$.}

\vspace{0,2cm}

Ideally, one would like then to interpret the different GIT quotients obtained (as $d$ decreases with respect to $g$) as first steps in a suitable ``Hassett-Keel" program for  Caporaso's 
compactified universal Jacobian $\ov J_{d,g}$ (see also Question \ref{QueB} and the discussion following it).  We hope to come back to this project in a future work.

In order to describe our results, we need to introduce some notation.

\section{Setup}

We work over an algebraically closed field $k$ (of arbitrary characteristic).
For an integer  $g\geq 2$ and any natural number  $d$, denote by $\Hilb_d$ the Hilbert scheme of curves
of degree $d$ and arithmetic genus $g$ in $\P^{d-g}=\P(V)$; denote by $\Chow_d$ the Chow scheme of $1$-cycles of degree $d$ in $\P^{d-g}$ and by
$$\Ch:\Hilb_d\to \Chow_d$$
the surjective map sending a one dimensional subscheme $[X\subset \P^{d-g}]\in \Hilb_d$ to its $1$-cycle. The linear algebraic group $\SL(V)\cong\SL_{d-g+1}$ acts naturally on
$\Hilb_d$ and $\Chow_d$ in such a way that $\Ch$ is an equivariant map. The action of $\SL(V)$ on $\Hilb_d$ and $\Chow_d$ can be naturally linearized as follows.

For any $m\gg0$, setting $P(m):=md+1-g$, the Hilbert scheme $\Hilb_d$ admits a $\SL(V)$-equivariant embedding
\begin{equation}\label{E:emb-Hilb}
j_m: \Hilb_d  \inj  \Gr(P(m),\Sym^m V^{\vee})  \inj \P\left(\bigwedge^{P(m)} \Sym^m V^{\vee}\right):=\P ,
\end{equation}
where $\Gr(P(m),\Sym^m V^{\vee})$ is the Grassmannian variety parametrizing $P(m)$-dimensional quotients of
$\Sym^m V^{\vee}$, which embeds in $\P\left(\bigwedge^{P(m)} \Sym^m V^{\vee}\right)$ via the Pl\"ucker embedding.
Explicitly, the map $j_m$ sends $[X\subset \P^{d-g}]\in \Hilb_d$ into
$$j_m([X\subset \P^{d-g}]):= \left[ \Sym^mV^{\vee} \surj H^0(X,\cO_X(m)) \right]
\in \Gr(P(m),\Sym^m V^{\vee})  \inj \P(\bigwedge^{P(m)} \Sym^m V^{\vee}).$$
We refer to  \cite[Lect. 15]{MumCAS} for details.
%It is well known that, for $m\gg0$,
%(see \cite[Lemma 2.1]{MS})
% any element $[X\subset \P^{d-g}]\in \Hilb_d$ is such that
%\begin{itemize}
%\item $\OO_X(m)$ has no higher cohomology;
%\item The natural map
%$$
%\Sym^m V^{\vee} \rightarrow \Gamma(\OO_X(m))=H^0(X, \OO_X(m))
%$$
%is surjective.
%\end{itemize}
%Under these hypothesis, if we set $P(m):=md+1-g$, then the $m$-th Hilbert point of $[X\subset \P^{d-g}]\in \Hilb_d$ is defined to be
%$$[X\subset \P^{d-g}]_m := \left[ \Sym^mV^{\vee} \surj \Gamma(\cO_X(m)) \right] \in
%\Gr(P(m),\Sym^mV^{\vee}) \inj \P\left(\bigwedge^{P(m)} \Sym^m V^{\vee}\right),$$
%where $\Gr(P(m),\Sym^m V^{\vee})$ is the Grassmannian variety parametrizing $P(m)$-dimensional quotients of
%$\Sym^m V^{\vee}$, which lies naturally in $\P\left(\bigwedge^{P(m)} \Sym^m V^{\vee}\right)$ via the Pl\"ucker embedding.
%\noindent In this way, for any $m\gg0$, we get a closed $\SL(V)$-equivariant embedding (see \cite[Lect. 15]{MumCAS}):
%$$
%\begin{array}{rcl}
%j_m: \Hilb_d & \inj & \Gr(P(m),\Sym^m V^{\vee})  \inj \P(\bigwedge^{P(m)} \Sym^m V^{\vee}):=\P \\
%\ X & \mapsto & [X]_m.
%\end{array}
%$$
From the embedding \eqref{E:emb-Hilb},  we get an ample $\SL(V)$-linearized line bundle
$\Lambda_m:=j_m^*\OO_{\P}(1)$ (for $m\gg0$) and we denote by
$$\Hilb_d^{s,m}\subseteq \Hilb_d^{ps,m}\subseteq \Hilb_d^{ss,m} \subseteq \Hilb_d$$
the locus of points that are, respectively, stable, polystable and semistable   with respect to
$\Lambda_m$.
%If $X\in \Hilb_d^{s,m}$ (resp. $X\in \Hilb_d^{ss,m}$), we say that $X$ is {\em $m$-Hilbert stable} (resp. {\em semistable}).
It is well-known (see \cite[Sec. 3.6]{HH2}) that
%The ample cone of $\Hilb$ admits a finite decomposition into locally-closed cells, such that the
%stable and the semistable locus are constant for linearizations taken from a given
%cell \cite[Theorem 0.2.3(i)]{DolHu}.
$\Hilb^{s,m}$, $\Hilb_d^{ps,m}$ and $\Hilb^{ss,m}$ are constant for $m\gg 0$ and we set
$$\begin{sis}
& \Hilb_d^s:=\Hilb_d^{s,m} & \text{ for } m\gg 0,\\
& \Hilb_d^{ps}:=\Hilb_d^{ps,m} & \text{ for } m\gg 0, \\
& \Hilb_d^{ss}:=\Hilb_d^{ss,m} & \text{ for } m\gg 0.
\end{sis}$$
If $[X\subset \P^{d-g}]\in \Hilb_d^{s}$ (resp. $ \Hilb_d^{ps}$, resp. $ \Hilb_d^{ss}$), we say that  $[X\subset \P^{d-g}]$ is {\em  Hilbert stable} (resp. {\em  Hilbert polystable},
resp. {\em  Hilbert semistable}).

The Chow scheme $\Chow_d$ admits a  $\SL(V)$-equivariant embedding
\begin{equation}\label{E:emb-Chow}
\Chow_d\stackrel{i}{\inj}\P(\otimes^{2} \Sym^d V^{\vee}):=\P'
\end{equation}
obtained by sending a $1$-cycle $Z$ of degree $d$ in $\P^{d-g}$ into the hyperplane of $\otimes^{2} \Sym^d V^{\vee}$ generated by all the elements $F\otimes G\in \Sym^d V^{\vee}$
such that  $Z\cap \{F=0\}\cap \{G=0\}\neq \emptyset$. We refer to \cite[Lec. 16]{MumCAS} for details.
From the embedding \eqref{E:emb-Chow}, we get  an ample $\SL(V)$-linearized line bundle $\Lambda:=i^*\OO_{\P'}(1)$
 and we denote by
$$\Chow_d^s\subseteq  \Chow_d^{ps}\subseteq \Chow_d^{ss} \subseteq \Chow_d$$
the locus of points of $\Chow_d$ that are, respectively, stable, polystable or semistable with respect to $\Lambda$.
%under the $\SL(V)$-action.
%We call the points of these loci, respectively, {\em Chow stable, polystable or semistable}. By a slight abuse of notation,
We say that $[X\subset \P^{d-g}]\in \Hilb_d$
 is \emph{Chow semistable} (resp. \emph{Chow polystable}, resp. \emph{Chow stable}) if its image $\Ch([X\subset \P^{d-g}])$ belongs to $\Chow_d^{ss}$ (resp. $\Chow_d^{ps}$,
 resp. $\Chow_d^s$).
%There is an $\SL(V)$-equivariant cycle map (see \cite[\S 5.4]{GIT}):
%$$
%\begin{array}{rcl}
%\Ch: \Hilb_d & \ra & \Chow_d  \\
 %        X & \mapsto & \Ch(X).
%\end{array}$$
The relation between Hilbert (semi)stability and Chow (semi)stability is given by the following chain of open inclusions (see \cite[Prop. 3.13]{HH2})
%By applying functoriality of stability \cite[Theorem 2.1]{Rei}, we get the following
\begin{equation}\label{E:rela-stability}
\Ch^{-1}(\Chow_d^s)\subseteq \Hilb_d^s\subseteq \Hilb_d^{ss}\subseteq \Ch^{-1}(\Chow_d^{ss})\subseteq \Hilb_d.
\end{equation}
%Let $[X]\in \Hilb_d$.
%\begin{enumerate}
%\item If $\Ch(X)\in \Chow_d^s$ then $[X]\in \Hilb_d^s$;
%\item If $[X]\in \Hilb_d^{ss}$ then $\Ch(X)\in \Chow_d^{ss}$.
%\end{enumerate}
In particular, there is a natural morphism of GIT-quotients
$$\Hilb_d^{ss}/\!\!/SL(V)\to \Chow_d^{ss}/\!\!/SL(V).$$
Note  that, in general, there is  no obvious relation between $\Hilb_d^{ps}$ and $\Ch^{-1}(\Chow_d^{ps})$: for example, according to \cite[Prop. 11.6 and Prop. 11.8]{HH2}, 
there are $2$-canonical curves that are Hilbert polystable but not Chow polystable and conversely.

We can now reformulate Problem(I) in the following form.

\vspace{0,2cm}

\textbf{Problem(II):}
\emph{Describe the points $[X\subset \P^{d-g}]\in \Hilb_d$ that are Hilbert or Chow (semi, poly) stable,
as $d$ decreases with respect to $g$.}

\vspace{0,2cm}

The aim of this note is to announce some partial results on the above Problem(II). Full details will appear in \cite{BMV}.

%We aim at giving a complete characterization of the GIT (semi-,poly-)stable points $[X\subset \P^{d-g}]\in \Hilb_d$ or of its image $\Ch([X\subset \P^{d-g}])\in \Chow_d$.

\section{Results}

%The aim of this research announcement is to state the results that we will prove

Our partial answer to the above Problem(II) will require some conditions on the singularities of $X$ and some conditions on the multidegree of the line bundle $\OO_X(1)$.
Let us introduce the relevant definitions.

%In order to describe our results, we need first to introduce some definitions.

%Let us list our results. First, let us first introduce some definitions.

\begin{defi}\label{D:quasi}
\noindent
\begin{enumerate}[(i)]
\item \label{quasi1} A curve $X$ is said to be \emph{quasi-stable}  if it is obtained from a stable curve $Y$ by ``blowing up" some of its nodes, i.e. by taking the partial normalization
of $Y$ at some of its nodes and inserting a $\P^1$ connecting the two branches of each node.
\item \label{quasi2} A curve $X$ is said to be \emph{quasi-p-stable}  if it is obtained from a p-stable curve $Y$ by ``blowing up" some of its nodes (as before) and ``blowing up" some of its cusps,
i.e. by taking the partial normalization  of $Y$ at some of its cusps and inserting a $\P^1$ tangent to the branch point of each cusp.
\end{enumerate}
Given a quasi-stable or a quasi-p-stable curve $X$, we call the $\P^1$'s inserted by blowing up nodes or cusps of $Y$ the \emph{exceptional components}, and we denote by $\exc\subset X$
the union of all of them.
\end{defi}

\begin{defi}
Let $X$ be a quasi-stable or a quasi-p-stable curve of genus $g\geq 2$ and let $L$ be a line bundle on $X$ of degree $d$. We say that:
\begin{enumerate}[(i)]
\item \label{bal1} $L$ is \emph{balanced} if for each subcurve $Z\subset X$ the following inequality
(called the basic inequality) is satisfied
\begin{equation*}
\left| \deg_ZL - \frac{d}{2g-2}\deg_Z(\omega_X)\right|\leq \frac{|Z\cap Z^c|}{2},\tag{*}
\end{equation*}
where $|Z\cap Z^c|$ denotes the length of the $0$-dimensional subscheme of $X$ obtained as the scheme-theoretic intersection of $Z$ with the complementary
subcurve $Z^c:=\ov{X\setminus Z}$.
\item \label{bal2} $L$ is  \emph{properly balanced}Ê if  $L$ is balanced and the degree of $L$ on each exceptional component of $X$ is $1$.
\item \label{bal3} $L$ is  \emph{strictly balanced} if $L$ is properly balanced and the basic inequality (*) is strict except possibly  for the subcurves $Z$ such that $Z\cap Z^c\subset \exc$.
\item \label{bal4} $L$ is  \emph{stably balanced} if $L$ is properly balanced and the basic inequality (*) is strict except possibly  for the subcurves $Z$ such that $Z$ or $Z^c$ is entirely contained
in $\exc$.
\end{enumerate}
\end{defi}

The inequality (*) first appeared in the work of Mumford \cite[Prop. 5.5]{Mum} and Gieseker \cite[Prop. 1.0.11]{Gie}. See also \cite[Sec. 3.1]{Cap}, where (*) 
is called the ``Basic Inequality".

We are now ready to state the main results of \cite{BMV}.
Our first result deals with high values of the degree $d$.

\begin{theoremalpha} \label{T:MainThmA}
Consider a point $[X\subset \P^{d-g}]\in \Hilb_d$ with $d>4(2g-2)$ and $g\geq 2$; assume moreover that $X$ is connected.
Then the following conditions are equivalent:
\begin{enumerate}[(i)]
\item $[X\subset \P^{d-g}] $ is Hilbert semistable (resp. polystable, resp. stable);
\item $[X\subset \P^{d-g}]$ is Chow semistable (resp. polystable, resp. stable);
\item $X$ is quasi-stable and $\OO_X(1)$ is balanced (resp. strictly balanced, resp. stably balanced).
\end{enumerate}
In each of the above cases, $X\subset \P^{d-g}$ is non-degenerate and linearly normal, and $\OO_X(1)$ is  non-special.

Moreover, the Hilbert or Chow GIT quotient is geometric (i.e. all the Hilbert or Chow semistable points are stable) if and only if $\gcd(2g-2, d-g+1)=1$.

\end{theoremalpha}

The above Theorem was proved by Caporaso in \cite{Cap} for $d\geq 10(2g-2)$ and for Hilbert (semi-, poly-)stability.
We remark that the hypothesis $d>4(2g-2)$ in the above Theorem \ref{T:MainThmA} is sharp: in \cite{HMo} it is  proved that a
$4$-canonical p-stable curve (which in particular can have cusps) is Hilbert stable  while a $4$-canonical stable curve with an
elliptic tail is not Hilbert semistable.

We then investigate what happens if $d\leq 4(2g-2)$ and we get a complete answer in the case $2(2g-2)<d<\frac{7}{2}(2g-2)$ and $g\geq 3$.

\begin{theoremalpha} \label{T:MainThmB}
Consider a point $[X\subset \P^{d-g}]\in \Hilb_d$ with $2(2g-2)<d<\frac{7}{2}(2g-2)$ and $g\geq 3$; assume moreover that $X$ is connected.
Then the following conditions are equivalent:
\begin{enumerate}[(i)]
\item $[X\subset \P^{d-g}] $ is Hilbert semistable (resp. polystable, resp. stable);
\item $[X\subset \P^{d-g}]$ is Chow semistable (resp. polystable, resp. stable);
\item $X$ is quasi-p-stable and $\OO_X(1)$ is balanced (resp. strictly balanced, resp. stably balanced).
\end{enumerate}
In each of the above cases, $X\subset \P^{d-g}$ is non-degenerate and linearly normal, and $\OO_X(1)$ is  non-special.

Moreover, the Hilbert or Chow GIT quotient is geometric (i.e. all the Hilbert or Chow semistable points are stable) if and only if $\gcd(2g-2, d-g+1)=1$.

\end{theoremalpha}

We note that the conditions on the degree $d$ and the genus $g$ in the above Theorem \ref{T:MainThmB} are sharp.
%lower bound $2(2g-2)<d$ and the upper bound $d<\frac{7}{2}(2g-2)$ on the degree is the sharpest condition under which the above Theorem holds true:
Indeed, if $d=2(2g-2)$ then it follows from  \cite[Thm. 2.14]{HH2} that there are $2$-canonical Hilbert stable curves having arbitrary tacnodes and not only tacnodes obtained by 
blowing up a cusp as in Definition \ref{D:quasi}\eqref{quasi2}.
%that a $2$-canonical h-stable curve in the sense of \cite[Def. 2.5, Def. 2.6]{HH2} 
%(which in particular can have arbitrary tacnodes and not only tacnodes obtained by blowing up a cusp as in Definition \ref{D:quasi}\eqref{quasi2})
%is Hilbert stable.
On the other hand, if $d=\frac{7}{2}(2g-2)$ (resp. $d>\frac{7}{2}(2g-2)$) then it follows from \cite[Prop. 1.0.6, Case 2]{Gie} that a point $[X\subset \P^{d-g}]\in \Hilb_d$
such that $X$ is a quasi-p-stable but not p-stable curve is not Chow stable (resp. Chow semistable)
\footnote{We thank Fabio Felici for pointing out to us the relevance of \cite[Prop. 1.0.6, Case 2]{Gie}.}.
Finally, if $g=3$ then Hyeon-Lee proved in \cite{HL} that a $3$-canonical irreducible p-stable curve  with one cusp is not Hilbert polystable.
 
%(while it is Hilbert semistable), which shows that the description of  Hilbert stable and Hilbert polystable points given in Theorem \ref{T:MainThmB} is false in this case. 
%Probably the description of Hilbert semistable points
%given in Theorem \ref{T:MainThmB} is still true for $g=2$; however for simplicity we restrict in this paper
 %to the case $g\geq 3$ whenever dealing with quasi-p-stable curves.

%The assumption  $g\geq 3$ allows us to avoid the subtleties that arise  while working with p-stable curves in genus $g=2$ (see \cite{HL}).

As an application of Theorem \ref{T:MainThmB}, we get a new compactification of the universal Jacobian    $J_{d,g}$ over the moduli space of p-stable curves of genus $g$.
%, which is the moduli scheme parametrizing
%pairs $(C,L)$ where $C$ is a smooth curve of genus $g$ and $L$ is a line bundle on $C$ of degree $d$.
To this aim, consider the category fibered in groupoids $\JJps$ over the category of schemes, whose fiber over a scheme $S$ is the groupoid of
families of quasi-p-stable curves over $S$ endowed with a line bundle whose restriction to the geometric fibers is properly balanced.

\begin{theoremalpha}\label{T:New-Comp-Jac}
Let $g\geq 3$ and $d\in \Z$.
\begin{enumerate}
\item $\JJps$ is a smooth,  irreducible  Artin stack of finite type over $k$ and of dimension $4g-4$. Moreover $\JJps$ is universally closed and
weakly separated (in the sense of \cite{ASvdW}).
%The category fibered in groupoids $\Pst_{d,g}$ representing families of quasi-p-stable curves endowed with a line bundle whose restriction to the geometric fibers is
%properly balanced is a smooth,  irreducible and universally closed Artin stack of dimension $4g-4$.

\item $\JJps$  admits an adequate moduli space $\Pdgps$ (in the sense of  \cite{alp}),
which is a normal irreducible projective variety  of dimension $4g-3$ containing $J_{d,g}$ as an open subvariety.
Moreover, if ${\rm char}(k)=0$, then $\Pdgps$ has rational singularities, hence it is Cohen-Macauly.

\item There exists a commutative digram
$$\xymatrix{
\JJps \ar[r] \ar[d]_{\Psi^{\rm ps}} & \Pdgps \ar[d]^{\Phi^{\rm ps}}\\
\ov{\mathcal M}_g^{\rm ps} \ar[r] & \Mgps
}$$
where $\Psi^{\rm ps}$ is surjective, universally closed and weakly separated (in the sense of \cite{ASvdW}) and $\Phi^{\rm ps}$ is surjective and projective with equidimensional
 fibers of dimension $g$.

%\item The following conditions are equivalent:
%\begin{enumerate}
%\item $\JJps$ is a Deligne-Mumford stack;
%\item $\JJps$ is proper;
%\item $\Psi^{\rm ps}$ is representable;
%\item $\Psi^{\rm ps}$ is proper;
%\item $ \JJps\to \Pdgps$ is a coarse moduli space;
%\item $\gcd(2g-2,d-g+1)=1$.
%\end{enumerate}

\item \label{T:New-Comp-Jac4} If ${\rm char}(k)=0$ or ${\rm char}(k)=p>0$ is bigger than the order of the automorphism group of any p-stable curve of genus $g$, then
for any $X\in \Mgps$, the fiber ${(\Phi^{\rm ps})}^{-1}(X)$  is isomorphic to  $\ov{\Jac_d}(X)/\Aut(X)$, where $\ov{\Jac_d}(X)$ is the Simpson's compactified Jacobian of $X$
parametrizing $S$-equivalence classes of rank $1$, torsion-free sheaves on $X$  that are slope-semistable with respect to $\omega_X$.

\item  If $2(2g-2)< d <\frac{7}{2}(2g-2)$ then $\JJps\cong [H_d/GL(r+1)]$ and $\Pdgps\cong H_d/GL(r+1)$, where $H_d\subset \Hilb_d$ is the open subset consisting of points
$[X\subset \P^{d-g}]\in \Hilb_d$ such that $X$ is connected and $[X\subset \P^{d-g}]$ is Hilbert semistable (or equivalently, Chow semistable).

\end{enumerate}
\end{theoremalpha}

The proof of the above Theorems will appear in \cite{BMV}.

\section{Open questions}

The above results leave unsolved some natural questions for further investigation, that we discuss briefly here.

The first question is of course the following

\begin{questionalpha}\label{QueA}
Describe the Hilbert and Chow (semi-,poly-)stable points of $\Hilb_d$  in the case where $\frac{7}{2}(2g-2)\leq d\leq 4(2g-2)$.
\end{questionalpha}
Indeed, as observed before, Theorems \ref{T:MainThmA} and \ref{T:MainThmB} do not hold for these values of $d$. In \cite[Thm. 5.1]{BMV}, we proved the following partial result:
if $[X\subset \P^{d-g}]$ is Chow semistable with $X$  connected then  $\OO_X(1)$ is balanced and $X$ is a reduced curve whose singularities are at most 
nodes, cusps and tacnodes at which two components of $X$ meet, one of which is a line of $\P^{d-g}$. Further progresses have been made by Fabio Felici in \cite{Fel}.

%one would expect that the GIT analysis gets harder as $d$ gets smaller. In particular, one would expect that the GIT analysis for
%$\frac{7}{2}(2g-2)\leq d\leq 4(2g-2)$ would be easier than for $2(2g-2)<d<\frac{7}{2}(2g-2)$, in which case we got a complete answer in Theorem \ref{T:MainThmB}.
%The problem however is that in the range $\frac{7}{2}(2g-2)\leq d< 4(2g-2)$ we have not been able to prove that if $[X\subset \P^{d-g}]\in \Hilb_d$ (with $X$ connected)
%is Chow semistable then $X$ does not have elliptic tails. Note however that this is certainly false for $d=4(2g-2)$ (see \cite{BMV} for details).

By analogy with the contraction map $T:\Mgps\to \Mgb$ constructed by Hassett-Hyeon in \cite{HH1}, the following problem arises very naturally.

\begin{questionalpha}\label{QueB}
Construct a map $\w{T}:\Jst \to \Jps$ fitting into the following commutative diagram
$$\xymatrix{
\Jst \ar[r]^{\w{T}} \ar[d]_{\Phi^s} & \Jps \ar[d]^{\Phi^{ps}} \\
\Mgb \ar[r]^TÊ& \Mgps,
}$$
where $\Jst\xrightarrow{\Phi^s} \Mgb$ is  Caporaso's compactification of the universal Jacobian over $\Mgb$.
%where $T$ is the map constructed by Hasset-Hyeon in \cite{HH1} (see also Fact \ref{F:mod-curves}\eqref{F:mod-curves2}).
\end{questionalpha}
More generally, one would like to set up a ``Hassett-Keel" program for
%the Caporaso's universal Jacobian
$\Jst$ and give an interpretation of the above map $\w{T}$ as the first step in this program.

Finally, the reader may have noticed that in Theorems \ref{T:MainThmA} and \ref{T:MainThmB} we have characterized points  $[X\subset \P^{d-g}]\in \Hilb_d$
that are Hilbert or Chow (semi,poly)stable under the assumption that $X$ is connected. Indeed, it is easy to prove that, at least for $d>2(2g-2)$,
the condition of being connected is both open and closed inside the Hilbert  or Chow semistable locus.
However, we do not know an answer to the following

\begin{questionalpha}\label{QueC}
Are there connected components inside the Hilbert or Chow GIT semistable locus made entirely of non-connected curves?
\end{questionalpha}

\end{document}